\documentclass[secthm,seceqn,amsthm,ussrhead,12pt]{amsart}
\usepackage[utf8]{inputenc}
\usepackage[english]{babel}

\usepackage{amssymb,amsmath,amsthm,amsfonts,xcolor,enumerate,hyperref,comment,longtable,cleveref}

\usepackage{times}
\usepackage{cite}
\usepackage{pdflscape}
\usepackage{ulem}
\usepackage[mathcal]{euscript}
\usepackage{tikz}
\usepackage{hyperref}
\usepackage{cancel}
\usepackage{stmaryrd}

\usetikzlibrary{arrows}

\newtheorem{theorem}{Theorem}
\newtheorem{lemma}[theorem]{Lemma}

\theoremstyle{definition}
\newtheorem{definition}[theorem]{Definition}

\theoremstyle{remark}
\newtheorem{remark}[theorem]{Remark}

\newenvironment{Proof}[1][Proof.]{\begin{trivlist}
\item[\hskip \labelsep {\bfseries #1}]}{\flushright
$\Box$\end{trivlist}}

\usepackage{stmaryrd}
\usepackage{xcolor}

\newcommand{\D}{\mathfrak{D}}

\newcommand{\la}{\langle}
\newcommand{\ra}{\rangle}

\newcommand{\nb}[1]{\nabla_{#1}}

\newcommand{\0}{\theta}
\newcommand{\af}{\alpha}

\begin{document}
\noindent{\Large  
The variety of dual mock-Lie algebras}\footnote{
The first part of this work is supported by the Russian Science Foundation under grant 19-71-10016. 
The second part of this work was supported by UAEU UPAR (9) 2017 Grant G00002599 (Mohamed A. Salim).
The authors thank  Prof. Dr. Pasha Zusmanovich for   discussions about dual mock-Lie algebras.

} \footnote{Corresponding Author: Ivan Kaygorodov }

   \

   {\bf    Luisa M. Camacho$^{a}$, Ivan Kaygorodov$^{b}$, Victor Lopatkin$^{c}$ \& Mohamed A. Salim$^{d}$   \\

    \medskip
 
    \medskip
}

{\tiny

$^{a}$ Dpto. Matem\'{a}tica Aplicada I, Universidad de Sevilla, Avda. Reina Mercedes, 41012 Sevilla, Spain

$^{b}$ CMCC, Universidade Federal do ABC, Santo Andr\'e, Brazil

$^{c}$ Saint Petersburg State University, Russia

$^{d}$ United Arab Emirates University, Al Ain, United Arab Emirates

\

\smallskip

  \medskip

   E-mail addresses: 
   
\smallskip

Luisa M. Camacho (lcamacho@us.es)

\smallskip
    
    Ivan Kaygorodov (kaygorodov.ivan@gmail.com) 

 \smallskip
 
    Victor Lopatkin (wickktor@gmail.com)

 \smallskip    
    Mohamed A. Salim (msalim@uaeu.ac.ae)

}

\ 

\ 

  \medskip

\ 

\noindent {\bf Abstract.}
{\it We classify all complex $7$- and $8$-dimensional dual mock-Lie algebras by algebraic and geometric way. Also we find all non-trivial complex $9$-dimensional dual mock-Lie algebras.}

\ 


\ 


\section*{Introduction}

There are many results related to the algebraic and geometric 
classification
of low-dimensional algebras in the varieties of Jordan, Lie, Leibniz and 
Zinbiel algebras;
for algebraic classifications  see, for example, \cite{ack,  degr3, usefi1, degr2, degr1, gkks, gkk, hac18, kkk18,  kv16};
for geometric classifications and descriptions of degenerations see, for example, 
\cite{ack, ale, ale2, BC99, gkks, gkk19, GRH, GRH2,  ikv17,  kkk18, kpv, kv16,kv17, S90}.
Here we give the algebraic and geometric classification of  
dual mock-Lie algebras of small dimensions. 

A while ago, a new class of algebras emerged in the literature -- the so-called mock-Lie algebras. These are commutative algebras satisfying the Jacobi identity. These algebras are locally nilpotent, so there are no nontrivial simple objects. Nevertheless, they seem to have an interesting structure theory which gives rise to interesting questions. And, after all, it is always curios to play with a classical notion by modifying it here and there and see what will happen – in this case, to replace in Lie algebras anti-commutativity by commutativity.

In \cite{hjs} (see also \cite{pasha}) that a finite-dimensional mock-Lie algebra does not admit a finite-dimensional faithful representation.
This class of algebras appeared in the literature under different names, reflecting, perhaps, the fact that it was considered from different viewpoints by different communities, sometimes not aware of each
other’s results. Apparently, for the first time these algebras appeared in \cite{zhev}, where an example of infinite-dimensional solvable but not nilpotent mock-Lie algebra was given (reproduced in \cite[ \S4.1, Example 1]{zsss}); further examples can be found in \cite[\S 4.1, Example 2 and \S 5.4, Exercise 4]{zsss} and \cite[\S 2.5]{w}. In this and other Jordan-algebraic literature these algebras are called just “Jordan algebras of
nil index 3”. In \cite{ok} they are called “Lie-Jordan algebras” (superalgebras are also considered there),
and, finally, in the recent papers \cite{bf} and \cite{am} the term “Jacobi–Jordan algebras” was used. The term
“mock-Lie” comes from \cite[\S 3.9]{gk}, where the corresponding operad appears in the list of quadratic
cyclic operads with one generator. Despite that Getzler and Kapranov downplayed this class of algebras
by (unjustly, in our opinion) calling them “pathological” (basing on the fact that the mock-Lie operad is
not Koszul), we prefer to stick to the term “mock-Lie”. These algebras live a dual life: as members of a
very particular class of Jordan algebras, and as strange cousins of Lie algebras

An entertaining fact (though not related to what follows): algebras over the operad Koszul dual to the mock-Lie operad can be characterized in three equivalent ways:
\begin{itemize}
\item anticommutative antiassociative algebras;
\item anticommutative $2$-Engel algebras;
\item anticommutative alternative algebras.
\end{itemize}
Here by antiassociative algebras we mean, following \cite{ok} and \cite{mr}, algebras satisfying the identity $(xy)z = -x(yz),$ and the $2$-Engel identity is $(xy)y = 0.$ Another entertaining fact (noted, for example, in \cite{ok}) is that mock-Lie algebras can be produced from antiassociative algebras the same way as they are produced from associative ones.

The algebraic classification of nilpotent algebras will be achieved by the calculation of central extensions of algebras from the same variety which have a smaller dimension.
Central extensions of algebras from various varieties were studied, for example, in \cite{ss78, zusmanovich, omirov}.
Skjelbred and Sund \cite{ss78} used central extensions of Lie algebras to classify nilpotent Lie algebras.
Using the same method,  
all non-Lie central extensions of  all $4$-dimensional Malcev algebras \cite{hac16},
all non-associative Jordan central extensions of all $3$-dimensional Jordan algebras,
all anticommutative central extensions of all  $3$-dimensional anticommutative algebras,
all central extensions of $2$-dimensional algebras,
and some others were described.
One can also look at the classification of
$3$-dimensional nilpotent algebras \cite{fkkv},
$4$-dimensional nilpotent associative algebras \cite{degr1},
$4$-dimensional nilpotent Novikov algebras \cite{kkk18},
$4$-dimensional nilpotent bicommutative algebras,
$4$-dimensional nilpotent commutative algebras in \cite{fkkv},
$5$-dimensional nilpotent restricted Lie agebras \cite{usefi1},
$5$-dimensional nilpotent Jordan algebras \cite{ha16},
$5$-dimensional nilpotent anticommutative algebras \cite{fkkv},
$6$-dimensional nilpotent Lie algebras \cite{degr3, degr2},
$6$-dimensional nilpotent Malcev algebras \cite{hac18},
$6$-dimensional nilpotent Tortkara algebras \cite{gkk,gkks},
$6$-dimensional nilpotent binary Lie algebras \cite{ack}.

Degenerations of algebras is an interesting subject, which has been studied in various papers.
In particular, there are many results concerning degenerations of algebras of small dimensions in a  variety defined by a set of identities.
One of important problems in this direction is a description of so-called rigid algebras. 
These algebras are of big interest, since the closures of their orbits under the action of the generalized linear group form irreducible components of the variety under consideration
(with respect to the Zariski topology). 
For example, rigid algebras in the varieties of
all  $4$-dimensional Leibniz algebras \cite{ikv17},
all  $4$-dimensional nilpotent Novikov algebras \cite{kkk18},
all  $4$-dimensional nilpotent bicommutative algebras,
all  $4$-dimensional nilpotent assosymmetric algebras,
all  $6$-dimensional nilpotent binary Lie algebras \cite{ack},
all  $6$-dimensional nilpotent Tortkara  algebras \cite{gkk19},
and in some other varieties were classified.
There are fewer works in which the full information about degenerations was given for some variety of algebras.
This problem was solved 
for $2$-dimensional pre-Lie algebras,  
for $2$-dimensional terminal algebras,
for $3$-dimensional Novikov algebras,  
for $3$-dimensional Jordan algebras,  
for $3$-dimensional Leibniz algebras, 
for $3$-dimensional anticommutative algebras,
for $3$-dimensional nilpotent algebras in \cite{fkkv},
for $4$-dimensional Lie algebras in \cite{BC99},
for $4$-dimensional Zinbiel  algebras,
for $4$-dimensional nilpotent Leibniz algebras,
for $4$-dimensional nilpotent commutative algebras in \cite{fkkv},
for $5$-dimensional nilpotent Tortkara algebras in \cite{gkks},
for $5$-dimensional nilpotent anticommutative algebras in \cite{fkkv},
for $6$-dimensional nilpotent Lie algebras in \cite{S90,GRH}, 
for $6$-dimensional nilpotent Malcev algebras in \cite{kpv}, 
for $2$-step nilpotent $7$-dimensional Lie algebras \cite{ale2}, 
and for all $2$-dimensional algebras in \cite{kv16}.

\section{The algebraic classification of  dual mock-Lie  algebras}


\subsection{The algebraic classification of [nilpotent] dual mock-Lie algebras}\label{algcl}
Let ${\bf A}$ and ${\bf V}$ be a dual mock-Lie algebra and a vector space and ${\rm Z}_{\D}^{2}\left( {\bf A},{\bf V}\right)$ denote the space of skew-symmetric  bilinear maps $\theta :{\bf A}\times 
{\bf A}\longrightarrow {\bf V}$
satisfying 
$ \theta(xy, z) =- \theta(x,yz).$
For $f\in{\rm Hom}({\bf A},{\bf V})$, we introduce $\delta f\in {\rm Z}_{\D}^{2}\left( {\bf A},{\bf V}\right)$ by the equality $\delta f\left( x,y\right) =f(xy)$ and
define ${\rm B}^{2}\left( {\bf A},{\bf V}\right) =\left\{\delta f \mid f\in {\rm Hom}\left( {\bf A},{\bf V}\right) \right\} $. One
can easily check that ${\rm B}^{2}({\bf A},{\bf V})$ is a linear subspace of ${\rm Z}_{\D}^{2}\left( {\bf A},{\bf V}\right)$.
Let us define $\rm {H}_{\D}^{2}\left( {\bf A},{\bf V}\right) $ as the quotient space ${\rm Z}_{\D}^{2}\left( {\bf A},{\bf V}\right) \big/{\rm B}^{2}\left( {\bf A},{\bf V}\right)$.
The equivalence class of $\theta\in {\rm Z}_{\D}^{2}\left( {\bf A%
},{\bf V}\right)$ in $\rm {H}_{\D}^{2}\left( {\bf A},{\bf V}\right)$ is denoted by $\left[ \theta \right]$. 

Suppose now that $\dim{\bf A}=m<n$ and $\dim{\bf V}=n-m$. For any dual mock-Lie 
bilinear map $\theta :{\bf A}\times {\bf A}\longrightarrow {\bf V%
}$, one can define on the space ${\bf A}_{\theta }:={\bf A}\oplus 
{\bf V}$ the dual mock-Lie bilinear product  $\left[ -,-\right] _{%
{\bf A}_{\theta }}$ by the equality $\left[ x+x^{\prime },y+y^{\prime }\right] _{%
{\bf A}_{\theta }}= xy  +\theta \left( x,y\right) $ for  
$x,y\in {\bf A},x^{\prime },y^{\prime }\in {\bf V}$. The algebra ${\bf A}_{\theta }$ is called an $(n-m)$-{\it %
dimensional central extension} of ${\bf A}$ by ${\bf V}$.
It is also clear that ${\bf A}_{\theta }$ is nilpotent if and only if so is ${\bf A}$.
The algebra ${\bf A}_{\theta }$ is dual mock-Lie  if
and only if ${\bf A}$ is dual mock-Lie  and $\theta$ is dual mock-Lie.

For a dual mock-Lie  bilinear form $\theta :{\bf A}\times {\bf A}\longrightarrow {\bf V}$, the space $\theta ^{\bot }=\left\{ x\in {\bf A}\mid \theta \left(
{\bf A},x\right) =0\right\} $ is called the {\it annihilator} of $\theta$.
For a dual mock-Lie  algebra ${\bf A}$, the ideal 
${\rm Ann}\left( {\bf A}\right) =\left\{ x\in {\bf A}\mid {\bf A}x =0\right\}$ is called the {\it annihilator} of ${\bf A}$.
One has
\begin{equation*}
{\rm Ann}\left( {\bf A}_{\theta }\right) =\left( \theta ^{\bot }\cap {\rm Ann}\left( 
{\bf A}\right) \right) \oplus {\bf V}.
\end{equation*}
Any $n$-dimensional  dual mock-Lie  algebra with non-trivial annihilator can be represented in
the form ${\bf A}_{\theta }$ for some $m$-dimensional  dual mock-Lie  algebra ${\bf A}$, an $(n-m)$-dimensional vector space ${\bf V}$ and $\theta \in {\rm Z}_{\D}^{2}\left( {\bf A},{\bf V}\right)$, where $m<n$ (see \cite[Lemma 5]{hac16}).
Moreover, there is a unique such representation with $m=n-\dim{\rm Ann}({\bf A})$. Note also that the last mentioned equality is equivalent to the condition  $\theta ^{\bot }\cap {\rm Ann}\left( 
{\bf A}\right)=0$. 

Let us pick some $\phi\in {\rm Aut}\left( {\bf A}\right)$, where ${\rm Aut}\left( {\bf A}\right)$ is the automorphism group of  ${\bf A}$.
For $\theta\in {\rm Z}_{\D}^{2}\left( {\bf A},{\bf V}\right)$, let us define $(\phi \theta) \left( x,y\right) =\theta \left( \phi \left( x\right)
,\phi \left( y\right) \right) $. Then we get an action of ${\rm Aut}\left( {\bf A}\right) $ on ${\rm Z}_{\D}^{2}\left( {\bf A},{\bf V}\right)$ that induces an action of the same group on $\rm {H}_{\D}^{2}\left( {\bf A},{\bf V}\right).$

\begin{definition}
Let ${\bf A}$ be an algebra and $I$ be a subspace of ${\rm Ann}({\bf A})$. If ${\bf A}={\bf A}_0 \oplus I$
then $I$ is called an {\it annihilator component} of ${\bf A}$.
\end{definition}

For a linear space $\bf U$, the {\it Grassmannian} $G_{s}\left( {\bf U}\right) $ is
the set of all $k$-dimensional linear subspaces of ${\bf U}$. For any $s\ge 1$, the action of ${\rm Aut}\left( {\bf A}\right)$ on $\rm {H}_{\D}^{2}\left( {\bf A},\mathbb{C}\right)$ induces 
an action of the same group on $G_{s}\left( \rm {H}_{\D}^{2}\left( {\bf A},\mathbb{C}\right) \right)$.
Let us define
$$
{\bf T}_{s}\left( {\bf A}\right) =\left\{ {\bf W}\in G_{s}\left( \rm {H}_{\D}^{2}\left( {\bf A},\mathbb{C}\right) \right)\left|\underset{[\theta]\in W}{\cap }\theta^{\bot }\cap {\rm Ann}\left( {\bf A}\right) =0\right.\right\}.
$$
Note that ${\bf T}_{s}\left( {\bf A}\right)$ is stable under the action of ${\rm Aut}\left( {\bf A}\right) $.

Let us fix a basis $e_{1},\ldots
,e_{s} $ of ${\bf V}$, and $\theta \in {\rm Z}_{\D}^{2}\left( {\bf A},{\bf V}\right) $. Then there are unique $\theta _{i}\in {\rm Z}_{\D}^{2}\left( {\bf A},\mathbb{C}\right)$ ($1\le i\le s$) such that $\theta \left( x,y\right) =\underset{i=1}{\overset{s}{%
\sum }}\theta _{i}\left( x,y\right) e_{i}$ for all $x,y\in{\bf A}$. Note that $\theta ^{\bot
}=\theta^{\bot} _{1}\cap \theta^{\bot} _{2}\cdots \cap \theta^{\bot} _{s}$ in this case.
If   $\theta ^{\bot
}\cap {\rm Ann}\left( {\bf A}\right) =0$, then by \cite[Lemma 13]{hac16} the algebra ${\bf A}_{\theta }$ has a nontrivial
annihilator component if and only if $\left[ \theta _{1}\right] ,\left[
\theta _{2}\right] ,\ldots ,\left[ \theta _{s}\right] $ are linearly
dependent in $\rm {H}_{\D}^{2}\left( {\bf A},\mathbb{C}\right)$. Thus, if $\theta ^{\bot
}\cap {\rm Ann}\left( {\bf A}\right) =0$ and the annihilator component of ${\bf A}_{\theta }$ is trivial, then $\left\langle \left[ \theta _{1}\right] , \ldots,%
\left[ \theta _{s}\right] \right\rangle$ is an element of ${\bf T}_{s}\left( {\bf A}\right)$.
Now, if $\vartheta\in {\rm Z}_{\D}^{2}\left( {\bf A},\bf{V}\right)$ is such that $\vartheta ^{\bot
}\cap {\rm Ann}\left( {\bf A}\right) =0$ and the annihilator component of ${\bf A}_{\vartheta }$ is trivial, then by \cite[Lemma 17]{hac16} one has ${\bf A}_{\vartheta }\cong{\bf A}_{\theta }$ if and only if
$\left\langle \left[ \theta _{1}\right] ,\left[ \theta _{2}%
\right] ,\ldots ,\left[ \theta _{s}\right] \right\rangle,
\left\langle \left[ \vartheta _{1}\right] ,\left[ \vartheta _{2}\right] ,\ldots,%
\left[ \vartheta _{s}\right] \right\rangle\in {\bf T}_{s}\left( {\bf A}\right)$ belong to the same orbit under the action of ${\rm Aut}\left( {\bf A}\right) $, where $%
\vartheta \left( x,y\right) =\underset{i=1}{\overset{s}{\sum }}\vartheta
_{i}\left( x,y\right) e_{i}$.

Hence, there is a one-to-one correspondence between the set of $%
{\rm Aut}\left( {\bf A}\right) $-orbits on ${\bf T}_{s}\left( {\bf A}%
\right) $ and the set of isomorphism classes of central extensions of $\bf{A}$ by $\bf{V}$ with $s$-dimensional annihilator and trivial annihilator component.
Consequently to construct all $n$-dimensional central extensions with $s$-dimensional annihilator and trivial annihilator component
of a given $(n-s)$-dimensional algebra ${\bf A}$ one has to describe ${\bf T}_{s}({\bf A})$, ${\rm Aut}({\bf A})$ and the action of ${\rm Aut}({\bf A})$ on ${\bf T}_{s}({\bf A})$ and then
for each orbit under the action of ${\rm Aut}({\bf A})$ on ${\bf T}_{s}({\bf A})$ pick a representative and construct the algebra corresponding to it. 

We will use the following auxiliary notation during the construction of central extensions.
Let ${\bf A}$ be an dual mock-Lie   algebra with  the basis $e_{1},e_{2},\ldots,e_{n}$. 
$\Delta_{ij}:{\bf A}\times {\bf A}\longrightarrow \mathbb{C}$ denotes the dual mock-Lie  bilinear form  defined by the equalities $\Delta _{ij}\left( e_{i},e_{j}\right)=-\Delta _{ij}\left( e_{j},e_{i}\right)=1$
and $\Delta _{ij}\left( e_{l},e_{m}\right) =0$ for $%
\left\{ l,m\right\} \neq \left\{ i,j\right\}$. In this case $\Delta_{ij}$ with $1\leq i < j\leq n $ form a basis of the space of dual mock-Lie  bilinear forms on $\bf{A}$.
We also denote by
$$\begin{array}{lll}
\D^i_{j}& \mbox{the }j\mbox{th }i\mbox{-dimensional dual mock-Lie algebra} \\
\end{array}$$

\subsection{The algebraic classification of low dimensional dual mock-Lie algebras}
Thanks to \cite{kkl19}, we have the classification of all $6$-dimensional nilpotent anticommutative algebras
and choosing only dual mock-Lie algebras from the list of algebras presented in \cite{kkl19}
we have the classification of all low dimensional dual mock-Lie algebras.
By the straightforward verification, it follows that only $\mathbb{M}_{01}$, $\mathbb{M}_{03}$, $\mathbb{M}_{04}$, $\mathbb{M}_{23}$, $\mathbb{M}_{24}$, and $\mathbb{M}_{26}$ satisfy antiassociativity low. We thus have the following table.

$$\begin{array}{lcllll}

\D^6_{01} 
&:& e_1e_2 = e_3    \\

\D^6_{02} 
&:& e_1e_2 = e_5 & e_3e_4 = e_5   \\

\D^6_{03} 
&:& e_1e_2 = e_4 & e_1e_3 = e_5   \\

\D^6_{04} 
&:& e_1e_3  =e_5 & e_2e_4=e_6  \\

\D^6_{05} 
&:& e_1e_2 =e_5 & e_1e_3 = e_6 & e_3e_4=e_5 \\

\D^6_{06} 
&:& e_1e_2=e_4 & e_1e_3=e_5 & e_2e_3=e_6

\end{array}$$

\subsection{The algebraic classification of $7$-dimensional dual mock-Lie algebras}
Thanks to \cite{ale2} we have the classification of all indecomposible $7$-dimensional $2$-step nilpotent dual mock-Lie algebras.

$$\begin{array}{lcllll}
\D^7_{07}  &:& e_1e_2=e_7 & e_3e_4=e_7 & e_5e_6=e_7  \\
\D^7_{08} &:& e_1e_2=e_6 & e_1e_4=e_7 & e_3e_5=e_7  \\
\D^7_{09} &:& e_1e_2=e_6 & e_1e_5=e_7 & e_3e_4=e_6 & e_2e_3=e_7 \\
\D^7_{10} &:& e_1e_2=e_5 & e_2e_3=e_6 & e_2e_4=e_7  \\
\D^7_{11} &:& e_1e_2=e_5 & e_2e_3=e_6 & e_3e_4=e_7  \\
\D^7_{12} &:& e_1e_2=e_5 & e_2e_3=e_6 & e_2e_4=e_7 & e_3e_4=e_5  \\
\D^7_{13} &:& e_1e_2=e_5 & e_1e_3=e_6 & e_2e_4=e_7 & e_3e_4=e_5  \\
\end{array}$$

The key tool in the classification of dual mock-Lie algebras will be the following obvious Lemma.

\begin{lemma}\label{splitext}
If the $i$-dimensional algebra $\D^i_{j}$ does not have  nontrivial dual mock-Lie central extension, 
then for every $k \in \mathbb N$ the $(i+k)$-dimensional algebra $\D^{i+k}_{j}$ does not have nontrivial dual mock-Lie central extensions.
\end{lemma}


Hence, for find non-$2$-step nilpotent $7$-dimensional dual mock-Lie algebras we need to calculate all non-split $2$-dimensional central extensions of all  $5$-dimensional dual mock-Lie algebras and 
all non-split $1$-dimensional central extensions of all $6$-dimensional dual mock-Lie algebras.
By some easy calculation, we have the cohomology spaces of these algebras.
\begin{equation*}
\begin{array}{|l|lll|l|} 
\hline
\mbox{$\D^5$}  & \multicolumn{3}{l|}{\mbox{ Multiplication table}} 
& \mbox{${\rm H}_{\D}^2({\D^5})$}  \\ 

\hline

{\D}_{01}^5& e_1e_2 = e_3 & & & \la [\Delta_{14}], [\Delta_{15}], [\Delta_{24}], [\Delta_{25}], [\Delta_{45}] \ra
  \\
\hline
{\D}_{02}^5 & e_1e_2=e_5 & e_3e_4=e_5 & &
\la[\Delta_{13}], [\Delta_{14}], [\Delta_{23}], [\Delta_{24}] \ra  \\

\hline

{\D}_{03}^5 & e_1e_2=e_4 & e_1e_3=e_5 & &
\la[\Delta_{23}]  \ra  \\

\hline

\multicolumn{5}{|c|}{}\\

\hline

\mbox{$\D^6$}  & \multicolumn{3}{l|}{\mbox{ Multiplication table}} 
& \mbox{${\rm H}_{\D}^2({\D^6})$}  \\






\hline
{\D}_{04}^6 & e_1e_3=e_5 & e_2e_4=e_6 & &
\la[\Delta_{12}],[\Delta_{14}], [\Delta_{23}], [\Delta_{34}] \ra  \\

\hline

{\D}_{05}^6 & e_1e_2=e_5 & e_1e_3=e_6 & e_3e_4=e_5 &
\la[\Delta_{14}], [\Delta_{23}], [\Delta_{24}] \ra  \\

\hline

{\D}_{06}^6 & e_1e_2=e_4 & e_1e_3=e_5 & e_2e_3=e_6 &
\la[\Delta_{16}]- [\Delta_{25}] + [\Delta_{34}] \ra  \\

\hline

\end{array}
\end{equation*}

Analizing the cohomology spaces of these algebras, we should conclude that only the algebra $\D^6_{06}$
has a non-split dual mock-Lie central extension.
Now, we have a new $7$-dimensional dual mock-Lie algebra

$$\begin{array}{lclllllll}

\D^7_{14} 
&:& e_1e_2=e_4 & e_1e_3=e_5 & e_1e_6=e_7 & e_2e_3=e_6 & e_2e_5=-e_7 & e_3e_4=e_7.

\end{array}$$

    \subsection{The algebraic classification of $8$-dimensional dual mock-Lie algebras}
It is easy to see that the algebra $\D^7_{14}$ has no non-trivial dual-mock-Lie central extentions. Hence, we will consider the cohomology space only for the following algebras.

    \begin{equation*}
    \begin{array}{|l|llll|l|} 
    \hline
    \mbox{$\D^7$}  & \multicolumn{4}{l|}{\mbox{ Multiplication table}}    & \mbox{${\rm H}_{\D}^2({\D^7})$}  \\ 
    \hline

    
    
 %

    


    {\D}_{06}^7 & e_1e_2=e_4& e_1e_3=e_5& e_2e_3=e_6& &
    \la[\Delta_{16}]- [\Delta_{25}] + [\Delta_{34}], [\Delta_{17}], [\Delta_{27}], [\Delta_{37}] \ra  \\
    \hline
    
    \D^7_{07}  & e_1e_2=e_7& e_3e_4=e_7& e_5e_6=e_7 && 
    \Big\la \begin{array}{l}[\Delta_{13}], [\Delta_{14}], [\Delta_{15}], [\Delta_{16}], [\Delta_{23}], [\Delta_{24}], \\

    [\Delta_{25}], [\Delta_{26}],[\Delta_{35}], [\Delta_{36}],[\Delta_{45}], [\Delta_{46}] \end{array} \Big\ra
    \\ 
    \hline
    
    \D^7_{08} & e_1e_2=e_6& e_1e_4=e_7&  e_3e_5=e_7& & \la [\Delta_{13}], [\Delta_{14}], [\Delta_{15}], [\Delta_{23}], [\Delta_{24}], [\Delta_{25}], [\Delta_{34}], [\Delta_{45}]
    \ra \\   
    \hline
    
    \D^7_{09} & e_1e_2=e_6& e_1e_5=e_7& e_3e_4=e_6& e_2e_3=e_7 & \la [\Delta_{13}], [\Delta_{14}], [\Delta_{24}], [\Delta_{25}], [\Delta_{35}], [\Delta_{45}]  \ra \\
    \hline
    
    \D^7_{10} & e_1e_2=e_5 & e_2e_3=e_6 & e_2e_4=e_7& & \la [\Delta_{13}], [\Delta_{14}], [\Delta_{34}]    \ra  \\
    \hline

    \D^7_{11} & e_1e_2=e_5& e_2e_3=e_6& e_3e_4=e_7 && \la [\Delta_{13}], [\Delta_{14}], [\Delta_{24}]    \ra  \\
    \hline
    
    \D^7_{12} & e_1e_2=e_5 & e_2e_3=e_6 & e_2e_4=e_7 & e_3e_4=e_5 & \la [\Delta_{13}], [\Delta_{14}]  \ra \\
    \hline
    
    \D^7_{13} & e_1e_2=e_5 & e_1e_3=e_6& e_2e_4=e_7& e_3e_4=e_5 & \la [\Delta_{14}], [\Delta_{23}]  \ra \\
    \hline
    

    \end{array}
    \end{equation*}
From here, only the algebra $\D^7_{06}$ maybe have a non-trivial dual mock-Lie central extension. We will find it.
The automorphism group $\mathrm{Aut}(\D_{06}^7)$ consists of invertible matrices of the form

\[
 \varphi = \left(\begin{matrix} 
  a & b & c & 0 & 0 & 0 & 0\\ 
  d & e & f & 0 & 0 & 0 & 0\\
  g & h & k & 0 & 0 & 0 & 0\\
  l & m & n & ae-db & af-dc & bf-ec & p \\
  q & r & s & ah-gb & ak-gc & bk-hc & i \\
  j & t & u & dh-ge & dk-gf & ek-hf & v \\
  w & x & y & 0 & 0 & 0 & z
 \end{matrix}\right).   
\]

Let us the notations
\[
 \nabla_1:=[\Delta_{16}]-[\Delta_{25}]+[\Delta_{34}], \quad \nabla_2:=[\Delta_{17}], \quad \nabla_3:=[\Delta_{27}], \quad \nabla_4:=[\Delta_{37}].
\]

Take $\theta=  \sum_{i=1}^4 \alpha_i\nabla_i \in \mathrm{H}^2_{\D}(\D_{06}^7,\mathbb{C})$. If $\varphi\in \mathrm{Aut}(\D_{06}^7)$, then

\[
 \varphi^T \left(\begin{matrix} 
 0 & 0 & 0 & 0 & 0& \alpha_1 & \alpha_2 \\
 0 & 0 & 0 & 0 & -\alpha_1& 0 & \alpha_3 \\
 0 & 0 & 0 & \alpha_1 & 0& 0 & \alpha_4 \\
 0 & 0 & -\alpha_1 & 0 & 0& 0 & 0 \\
 0 & \alpha_1 & 0 & 0 & 0& 0 & 0 \\
 -\alpha_1 & 0 & 0 & 0 & 0& 0 & 0 \\
 -\alpha_2 & -\alpha_3 & -\alpha_4 & 0 & 0& 0 & 0 \\
 \end{matrix}\right)\varphi =
 \left( \begin{matrix} 0 & \beta_1^* & \beta_2^*& 0 & 0 & \alpha_1^* & \alpha_2^* \\
 -\beta_1^* & 0 & \beta_3^* & 0 & -\alpha_1^* & 0 & \alpha_3^*\\
 -\beta_2^* & - \beta_3^* & 0 & \alpha_1^* & 0 & 0 & \alpha_4^* \\
 0 & 0 & -\alpha_1^* & 0&0&0&0 \\
 0 & \alpha_1^* & 0&0&0&0&0\\
 -\alpha_1^* & 0&0&0&0&0&0 \\
 - \alpha_2^* & -\alpha_3^* & - \alpha_4^* & 0&0&0&0
 \end{matrix}  \right)
\]

\begin{align*}
 & \alpha^*_{1} = -(c e g - b f g - c d h + a f h + b d k - a e k) \alpha_1,\\
 & \alpha_2^* = (-d i  + g p + a v)\alpha_1 + a z \alpha_2 + d z \alpha_3 + g z \alpha_4,\\
 & \alpha_3^* = (-e i  + h p  + b v ) \alpha_1 + b z \alpha_2 + e z \alpha_3 + h z \alpha_4, \\
 & \alpha_4^* = (-f i + k p + c v) \alpha_1 + c z \alpha_2 + f z \alpha_3 + k z \alpha_4.
\end{align*}


Hence, $\phi\langle\0\rangle=\langle\0^*\rangle,$ where $\0^*=\sum_{i=1}^4 \af_i^*  \nb i.$ 
We are interesting in elements with $\alpha_1 \neq 0$ and $(\alpha_2, \alpha_3, \alpha_4) \neq (0,0,0).$ 
Without loss of generality, we can suppose that $\alpha_4\neq 0.$
So, by choosing the following non-zero elements $d=h=a=e=k=1$ and 
\[ 
v=1-\dfrac{\alpha_2}{\alpha_4}, \
i =1+ \dfrac{\alpha_3}{\alpha_4},\
z= \dfrac{\alpha_1}{\alpha_4},\ 
c= -1 + \dfrac{1}{\alpha_1}
 \]
 we get the representative $\langle \nabla_1 + \nabla_4\rangle$. Now we have the new $8$-dimensional dual mock-Lie algebra $\D_{36}^8$ constructed from $\D_{06}^7$:
$$\begin{array}{lclllllll}
 \D_{36}^8&:& e_1e_2 = e_4& e_1e_3 = e_5& e_2e_3 = e_6& e_1e_6 = e_8& e_2e_5=-e_8& e_3e_4 = e_8& e_3e_7 = e_8.
\end{array}
$$

Thanks to \cite{ale} we have the list of all $8$-dimensional $2$-step nilpotent indecomposible Lie algebras:

$$\begin{array}{lclllllll}









\D_{15}^8 &:& e_1e_2 = e_4 & e_3e_2 = e_5 & e_6e_7=e_8 \\

\D_{16}^8 &:& e_1e_2 = e_5 & e_3e_4 = e_5 & e_6e_7=e_8 \\

\D_{17}^8 &:& e_1e_2=e_7 &  e_3e_4=e_8 & e_5e_6=e_7+e_8 \\

\D_{18}^8 &:& e_1e_2=e_7 & e_4e_5=e_7 & e_1e_3=e_8 & e_4e_6=e_8\\

\D_{19}^8 &:& e_1e_2=e_7 & e_4e_5=e_7 & e_3e_4=e_8 & e_5e_6=e_8\\

\D_{20}^8 &:& e_1e_2=e_7 & e_3e_4=e_7 & e_5e_6=e_7& e_4e_5=e_8\\

\D_{21}^8 &:& e_1e_2=e_7 & e_3e_4=e_7 & e_5e_6=e_7 & e_2e_3=e_8 & e_4e_5=e_8 \\

\D_{22}^8 &:& e_1e_2=e_6 & e_4e_5=e_6 & e_2e_3=e_7 & e_1e_3=e_8 \\

\D_{23}^8 &:& e_1e_2=e_6 & e_4e_5=e_6 & e_2e_3=e_7 & e_3e_4=e_8 \\

\D_{24}^8 &:& e_1e_2=e_6 & e_2e_3=e_7 & e_4e_5=e_7 & e_3e_4=e_8  \\

\D_{25}^8 &:& e_1e_2=e_6 & e_2e_3=e_7 & e_4e_5=e_7 & e_3e_4=e_8 &  e_5e_1=e_8 \\ 

\D_{26}^8 &:& e_1e_2=e_6 & e_1e_3=e_7 & e_1e_4=e_8 & e_2e_5=e_7 \\

\D_{27}^8 &:& e_1e_2=e_6 & e_1e_3=e_7 & e_1e_4=e_8 & e_2e_3=e_8 & e_4e_5=e_7  \\

\D_{28}^8 &:& e_1e_2=e_6 & e_1e_3=e_7 & e_1e_5=e_8 & e_2e_4=e_8 & e_3e_4=e_6\\

\D_{29}^8 &:& e_1e_2=e_6 & e_1e_3=e_7 & e_2e_3=e_8 & e_1e_4=e_8 & e_2e_5=e_7 \\

\D_{30}^8 &:& e_1e_2=e_6 & e_1e_3=e_7 & e_2e_3=e_8 & e_1e_4=e_8 & e_2e_5=e_7 & e_4e_5=e_6 \\ 

\D_{31}^8 &:& e_1e_2=e_6 & e_2e_3=e_7 & e_3e_4=e_7 & e_4e_5=e_8\\ 

\D_{32}^8 &:& e_1e_2=e_6 & e_2e_3=e_7 & e_3e_4=e_8 & e_4e_5=e_7 & e_5e_1=e_7\\ 

\D_{33}^8 &:& e_1e_2=e_5 & e_2e_3=e_6 & e_3e_4=e_7 & e_4e_1=e_8 \\

\D_{34}^8 &:& e_1e_2=e_5 & e_1e_3=e_6 & e_2e_3=e_7 & e_1e_4=e_8 \\

\D_{35}^8 &:& e_1e_2=e_5 & e_1e_3=e_6 & e_2e_4=e_6 & e_2e_3=e_7 & e_1e_4=e_8\\

\end{array}$$

\subsection{The algebraic classification of $9$-dimensional  dual mock-Lie algebras}
The description of $2$-step nilpotent Lie algebras is not finished now.
There is only some particular classification of these algebras \cite{ren}.
Here, we give the classification of all complex $9$-dimensional non-Lie dual mock-Lie algebras.
Analyzing the dimension of cohomology spaces of $i$-dimensional $2$-step nilpotent Lie algebras ($i=3,4,5,6,7)$, 
we conclude that only $\D^7_{06}$ maybe give some non-trivial $(9-i)$-dimensional dual mock-Lie cenrtal extensions.
Hence, we will calculate $2$-dimensional dual mock-Lie central extensions of $\D^7_{06}$ and $1$-dimensional dual mock-Lie extensions of $8$-dimensional $2$-step nilpotent Lie algebras.

\subsubsection{$2$-dimensional dual mock-Lie central extensions of $7$-dimensional $2$-step nilpotent Lie algebras}
Here we are considering $2$-dimensional dual mock-Lie central extensions of $\D^7_{06}.$
Consider the vector space generated by the following two cocycles 
$$\begin{array}{rcl}
\theta_1 &=& \alpha_1 \nabla_1+\alpha_2\nabla_2+\alpha_3\nabla_3+\alpha_4\nabla_4   \\
\theta_2 &=& \beta_2 \nabla_2+\beta_3\nabla_3+\beta_4\nabla_4.
\end{array}$$
It is easy to see, that we can suppose that $\alpha_1 \beta_2 \neq 0.$
Then  by choosing 
    the following nonzero elements 
\[ 
a = -\frac{\beta_3}{\alpha_1}, \ 
b = -\frac{\beta_4}{\beta_2}, \
c = \frac{1}{\beta_2}, \  
d = \frac{\beta_2}{\alpha_1}, \ 
h = 1, \ 
i = \frac{\alpha_3}{\alpha_1}, \
p = -\frac{\alpha_4}{\alpha_1}, \ 
v = -\frac{\alpha_2}{\alpha_1}, \ 
z = 1,
 \] 
we have the representative $\langle \nabla_1, \nabla_4 \rangle$ which gives the following $9$-dimensional algebra:

$$\begin{array}{lclllllll}
 \D_{37}^9 &:& e_1e_2 = e_4& e_1e_3 = e_5& e_2e_3 = e_6& e_1e_6 = e_8& e_2e_5=-e_8& e_3e_4 = e_8& e_3e_7 = -e_9.
\end{array}$$

\subsubsection{$1$-dimensional dual mock-Lie central extensions of $8$-dimensional $2$-step nilpotent Lie algebras}

By Lemma \ref{splitext} and \cite[Theorem 3.8, 3.9]{ale}, we have the following dual mock-Lie algebras have nontrivial dual mock-Lie extensions.

    \begin{longtable}{|l|l|} 
    \hline
    $\D^8$    
    &  ${\rm H}_{\D}^2({\D^8})$  \\ 
    
    \hline
    
    $\D^8_{06}$ &  $\la [\Delta_{16}] - [\Delta_{25}]+[\Delta_{34}], [\Delta_{17}],[\Delta_{18}], [\Delta_{27}], [\Delta_{28}], [\Delta_{37}], [\Delta_{38}], [\Delta_{78}]\ra$      \\
    \hline
    
    $ \D_{15}^8$ & $\la [\Delta_{13}], [\Delta_{16}], [\Delta_{17}], [\Delta_{26}], [\Delta_{27}], [\Delta_{36}], [\Delta_{37}]  \ra$ \\
    \hline
    
   $ \D_{16}^8 $  &  $ \la [\Delta_{13}], [\Delta_{14}], [\Delta_{16}], [\Delta_{17}], [\Delta_{23}], [\Delta_{24}], [\Delta_{26}], [\Delta_{27}], [\Delta_{36}], [\Delta_{37}], [\Delta_{46}], [\Delta_{47}] \ra$  \\
   \hline
   
   $ \D_{17}^8$  & $ \la [\Delta_{13}], [\Delta_{14}], [\Delta_{15}], [\Delta_{16}], [\Delta_{23}], [\Delta_{24}], [\Delta_{25}], [\Delta_{26}], [\Delta_{35}], [\Delta_{36}], [\Delta_{45}], [\Delta_{46}]  \ra$  \\
   \hline
   
  $ \D_{18}^8$  &  $ \la [\Delta_{14}], [\Delta_{15}], [\Delta_{16}], [\Delta_{23}], [\Delta_{24}], [\Delta_{25}], [\Delta_{26}], [\Delta_{34}], [\Delta_{35}], [\Delta_{36}], [\Delta_{56}] \ra$  \\
  \hline
       
    $ \D_{19}^8$  &  $ \la [\Delta_{13}],[\Delta_{14}], [\Delta_{15}], [\Delta_{16}], [\Delta_{23}], [\Delta_{24}], [\Delta_{25}], [\Delta_{26}], [\Delta_{35}], [\Delta_{36}], [\Delta_{46}]  \ra$  \\
 \hline
 
 $ \D_{20}^8$  &$   \la [\Delta_{13}], [\Delta_{14}], [\Delta_{15}], [\Delta_{16}], [\Delta_{23}], [\Delta_{24}], [\Delta_{25}], [\Delta_{13}], [\Delta_{26}], [\Delta_{35}], [\Delta_{36}], [\Delta_{46}] \ra $ \\
 \hline
  
 $ \D_{21}^8$  & $  \la [\Delta_{13}], [\Delta_{14}], [\Delta_{15}], [\Delta_{16}], [\Delta_{24}], [\Delta_{25}], [\Delta_{26}], [\Delta_{35}], [\Delta_{36}], [\Delta_{46}] \ra$  \\

 \hline

$ \D_{22}^8$  & $ \la [\Delta_{14}], [\Delta_{15}], [\Delta_{24}], [\Delta_{25}], [\Delta_{34}], [\Delta_{35}] \ra $ \\
\hline

$ \D_{23}^8$  & $  \la [\Delta_{13}], [\Delta_{14}], [\Delta_{15}], [\Delta_{24}], [\Delta_{25}], [\Delta_{35}] \ra$  \\
\hline

$ \D_{24}^8$  & $  \la [\Delta_{13}], [\Delta_{14}], [\Delta_{15}], [\Delta_{24}], [\Delta_{25}], [\Delta_{35}] \ra$  \\
\hline

$ \D_{25}^8$  & $  \la [\Delta_{13}], [\Delta_{14}], [\Delta_{24}], [\Delta_{25}], [\Delta_{35}] \ra $ \\
\hline

$ \D_{26}^8$  & $ \la [\Delta_{15}], [\Delta_{23}], [\Delta_{24}], [\Delta_{34}], [\Delta_{35}], [\Delta_{45}] \ra $ \\
\hline

$ \D_{27}^8 $ & $ \la [\Delta_{15}], [\Delta_{24}], [\Delta_{15}], [\Delta_{34}],[\Delta_{35}]  \ra $ \\
\hline

$ \D_{28}^8 $ & $  \la [\Delta_{14}], [\Delta_{23}], [\Delta_{25}], [\Delta_{35}], [\Delta_{45}] \ra$ \\
\hline

$ \D_{29}^8$  & $ \la [\Delta_{15}], [\Delta_{24}], [\Delta_{34}], [\Delta_{35}], [\Delta_{45}] \ra $ \\
\hline

$ \D_{30}^8$  & $ \la [\Delta_{15}], [\Delta_{24}], [\Delta_{34}], [\Delta_{35}]  \ra $  \\
\hline

$ \D_{31}^8 $ & $ \la [\Delta_{13}], [\Delta_{14}], [\Delta_{15}], [\Delta_{24}],[\Delta_{25}],[\Delta_{35}] \ra$  \\
\hline

$ \D_{32}^8 $ &   $ \la [\Delta_{13}], [\Delta_{14}], [\Delta_{24}], [\Delta_{25}], [\Delta_{35}]  \ra$ \\
\hline

$ \D_{33}^8 $ & $  \la [\Delta_{13}], [\Delta_{24}]  \ra $ \\
\hline

$ \D_{34}^8 $ & $  \la [\Delta_{15}], [\Delta_{24}], [\Delta_{25}], [\Delta_{34}], [\Delta_{35}], [\Delta_{45}]  \ra $ \\
\hline

$ \D_{35}^8 $ & $   \la [\Delta_{34}] \ra $ \\
\hline
\end{longtable} 
 From here, only the algebra $\D^8_{06}$ maybe have a non-trivial dual mock-Lie central extension. We will find it.   
The automorphism group $\mathrm{Aut}(\D_{06}^8)$ consists of invertible matrices of the form

\[
 \varphi = \left(\begin{matrix} 
  a & b & c & 0 & 0 & 0 & 0 & 0\\ 
  d & e & f & 0 & 0 & 0 & 0 & 0\\
  g & h & k & 0 & 0 & 0 & 0 & 0\\
  l & m & n & ae-db & af-dc & bf-ec & p_1 & p_2 \\
  q & r & s & ah-gb & ak-gc & bk-hc & i_1 & i_2 \\
  j & t & u & dh-ge & dk-gf & ek-hf & v_1 & v_2 \\
  w_1 & x_1 & y_1 & 0 & 0 & 0 & z_1 & z_2 \\
  w_2 & x_2 & y_2 & 0 & 0 & 0 & z_3 & z_4
 \end{matrix}\right).
\]

Let us use the notations
\[ \nabla_1 := [\Delta_{16}]-[\Delta_{25}]+[\Delta_{34}],\ 
\nabla_2 := [\Delta_{17}], \
\nabla_3 := [\Delta_{18}], \]
\[ \nabla_4 := [\Delta_{27}], \
\nabla_5 := [\Delta_{28}], \
\nabla_6 := [\Delta_{37}], \
\nabla_7 := [\Delta_{38}], \
\nabla_8 := [\Delta_{78}]. \] 

Take $\theta=  \sum_{i=1}^8 \alpha_i\nabla_i \in \mathrm{H}^2_{\D}(\D_{06}^8,\mathbb{C})$. If $\varphi\in \mathrm{Aut}(\D_{06}^8)$, then

\[
 \varphi^T \left(\begin{matrix} 
 0 & 0 & 0 & 0 & 0& \alpha_1 & \alpha_2 &\alpha_3 \\
 0 & 0 & 0 & 0 & -\alpha_1& 0 & \alpha_4 & \alpha_5 \\
 0 & 0 & 0 & \alpha_1 & 0& 0 & \alpha_6 & \alpha_7\\
 0 & 0 & -\alpha_1 & 0 & 0& 0 & 0 & 0\\
 0 & \alpha_1 & 0 & 0 & 0& 0 & 0 & 0 \\
 -\alpha_1  & 0 & 0 & 0 & 0& 0 & 0 & 0 \\
 -\alpha_2 & -\alpha_4 & -\alpha_6 & 0 & 0& 0 & 0 & \alpha_8 \\
-\alpha_3 & -\alpha_5 & -\alpha_7 & 0 & 0& 0 & -\alpha_8 &0\\
 \end{matrix}\right)\varphi =
 \left( \begin{matrix}
 0 & \beta_1^* & \beta_2^*& 0 & 0 & \alpha_1^* & \alpha_2^* & \alpha_3^* \\
 -\beta_1^* & 0 & \beta_3^* & 0 & -\alpha_1^* & 0 & \alpha_4^* & \alpha_5^*\\
 -\beta_2^* & - \beta_3^* & 0 & \alpha_1^* & 0 & 0 & \alpha_6^* & \alpha_7* \\
 0 & 0 & -\alpha_1^* & 0&0&0&0&0 \\
 0 & \alpha_1^* & 0&0&0&0&0&0\\
 -\alpha_1^* & 0&0&0&0&0&0&0 \\
 - \alpha_2^* & -\alpha_4^* & - \alpha_6^* & 0&0&0&0&\alpha^*_8\\
- \alpha_3^* & -\alpha_5^* & - \alpha^*_7 & 0&0&0&-\alpha_8^* & 0
 \end{matrix}  \right),
\]
where
\begin{eqnarray*}
 \alpha_1^* &= & -(c e g - b f g - c d h + a f h + b d k - a e k) \alpha_1, \\
 \alpha_2^* &=& (-d  i_1 + g p_1 + 
 a v_1)\alpha_1 + (a \alpha_2 + d \alpha_4 + g \alpha_6 -  w_2\alpha_8) z_1 + (a \alpha_3 + d \alpha_5 + g \alpha_7 +  w_1\alpha_8) z_3, \\
 \alpha_3^* &=& (-d i_2 + g p_2 + 
  a v_2)\alpha_1 + (a \alpha_2 + d \alpha_4 + g \alpha_6 -w_2\alpha_8) z_2 + (a \alpha_3 + d \alpha_5 + g \alpha_7 + w_1\alpha_8) z_4 \\
  \alpha_4^* &=& (-e i_1 + h p_1 + 
 b v_1)\alpha_1 + (b \alpha_2 + e \alpha_4 + h \alpha_6 -  x_2\alpha_8)z_1 + (b\alpha_3 + e \alpha_5 + h \alpha_7 +  x_1\alpha_8) z_3, \\
 \alpha_5^* &=& (e i_2 + h p_2 + 
  b v_2)\alpha_1 + (b \alpha_2 + e \alpha_4 + h \alpha_6 - x_2\alpha_8) z_2 + (b \alpha_3 + e\alpha_5 + h\alpha_7 +  x_1\alpha_8) z_4, \\
  \alpha_6^* &=& (-f i_1 + k p_1 + 
 c v_1) \alpha_1 + (c \alpha_2 + f \alpha_4 + k \alpha_6 -  y_1\alpha_8) z_1 + (c \alpha_3 + f \alpha_5 + k \alpha_7 + y_1\alpha_8) z_3,\\
 \alpha_7^* &=& (-f i_2 + k p_2 + 
  c v_2)\alpha_1 + (c \alpha_2 + f \alpha_4 + k\alpha_6 - y_1\alpha_8) z_2 + (c\alpha_3 + f\alpha_5 + k \alpha_7+  y_1\alpha_8) z_4,\\
  \alpha_8^* &=& (-z_2 z_3 + z_1 z_4)\alpha_8.
\end{eqnarray*}

Hence, $\phi\langle\0\rangle=\langle\0^*\rangle,$ where $\0^*=\sum_{i=1}^8 \af_i^*  \nb i.$
Here, we have the following situations:

\begin{enumerate}
    \item $\alpha_1, \alpha_8 \neq 0,$ then by choosing 
    the following nonzero elements 
\[ c = 1,\ e = 1, \ h = 1,\ k = 1, \ g =  - \frac{1}{\alpha_1}, \ 
z_1 = \frac{1}{\alpha_8}, \ 
z_2 = 2, \ z_4 = 1,\] 
\[ y_1 = \frac{-\alpha_3 +\alpha_5}{\alpha_8}, \  
x_2 = \frac{-\alpha_2 - \alpha_3 + \alpha_4  + \alpha_5}{\alpha_8},\ 
p_1 = -\frac{\alpha_2  + \alpha_3  - \alpha_5  +\alpha_6}{\alpha_1 \alpha_8}, \]
\[ p_2 = -\frac{ 2 \alpha_2 + 2 \alpha_3 - \alpha_5  + 2 \alpha_6  +\alpha_7}{\alpha_1}, \ 
w_1 = -\frac{\alpha_5 }{\alpha_1\alpha_8}, \
w_2 = \frac{\alpha_2  +\alpha_3 -\alpha_5}{ \alpha_1 \alpha_8},\]
we have the representative $\langle \nabla_1+\nabla_8 \rangle.$  
Now we have the new $9$-dimensional dual mock-Lie algebra:
\[
 \D_{38}^9: 
 e_1e_2 = e_4, e_1e_3 = e_5, e_2e_3 = e_6, e_1e_6 = e_9, e_2e_5=-e_9, e_3e_4 = e_9, e_7e_8 = e_9.
\]

\item $\alpha_1\neq 0, \alpha_8=0,$ then by choosing the following nonzero elements 

\[ a=\frac{1}{\alpha_1},\ e=1,\ k=1,\
v_1 = -\frac{\alpha_2}{\alpha_1}, \ 
v_2 = -\frac{\alpha_3}{\alpha_1}, \ 
p_1 = -\frac{\alpha_6}{\alpha_1}, \
p_2 = -\frac{\alpha_7}{\alpha_1}, \
i_1 = \frac{\alpha_4}{\alpha_1}, \
i_2 = \frac{\alpha_5}{\alpha_1}, \]
we have the representative $\langle \nabla_1 \rangle$ and it is a split algebra.

\item if $\alpha_1, \alpha_8=0,$ then we can suppose that $\alpha_7\neq 0$ and by choosing the following nonzero elements 

\[
a=1, \ e = 1,\  k = 1, \ f = 1, \ z_1=1, \ z_4=1,\
g = -\frac{\alpha_2}{\alpha_6}, \
h = -\frac{\alpha_4}{\alpha_6}, \ 
k = -\frac{\alpha_4}{\alpha_6},
\] 
then we have a representative from $\langle \nabla_2, \nabla_4 , \nabla_6 \rangle,$ which gives a split algebra.
\end{enumerate}

Summarizing, we have the following theorem

\begin{theorem}
Let $\D$ be a complex $9$-dimensional indecomposible non-Lie dual mock-Lie algebra,
then $\D$ is isomorphic to $\D^9_{37}$ or $\D^9_{38}.$
\end{theorem}
\section{Degenerations of   dual mock-Lie algebras }

\subsection{Degenerations of algebras}
Given an $n$-dimensional vector space ${\bf V}$, the set ${\rm Hom}({\bf V} \otimes {\bf V},{\bf V}) \cong {\bf V}^* \otimes {\bf V}^* \otimes {\bf V}$ 
is a vector space of dimension $n^3$. This space has a structure of the affine variety $\mathbb{C}^{n^3}.$ 
Indeed, let us fix a basis $e_1,\dots,e_n$ of ${\bf V}$. Then any $\mu\in {\rm Hom}({\bf V} \otimes {\bf V},{\bf V})$ is determined by $n^3$ structure constants $c_{i,j}^k\in\mathbb{C}$ such that
$\mu(e_i\otimes e_j)=\sum_{k=1}^nc_{i,j}^ke_k$. A subset of ${\rm Hom}({\bf V} \otimes {\bf V},{\bf V})$ is {\it Zariski-closed} if it can be defined by a set of polynomial equations in the variables $c_{i,j}^k$ ($1\le i,j,k\le n$).

Let $T$ be a set of polynomial identities.
All algebra structures on ${\bf V}$ satisfying polynomial identities from $T$ form a Zariski-closed subset of the variety ${\rm Hom}({\bf V} \otimes {\bf V},{\bf V})$. We denote this subset by $\mathbb{L}(T)$.
The general linear group ${\rm GL}({\bf V})$ acts on $\mathbb{L}(T)$ by conjugation:
$$ (g * \mu )(x\otimes y) = g\mu(g^{-1}x\otimes g^{-1}y)$$ 
for $x,y\in {\bf V}$, $\mu\in \mathbb{L}(T)\subset {\rm Hom}({\bf V} \otimes {\bf V},{\bf V})$ and $g\in {\rm GL}({\bf V})$.
Thus, $\mathbb{L}(T)$ is decomposed into ${\rm GL}({\bf V})$-orbits that correspond to the isomorphism classes of algebras. 
Let $O(\mu)$ denote the ${\rm GL}({\bf V})$-orbit of $\mu\in\mathbb{L}(T)$ and $\overline{O(\mu)}$ its Zariski closure.

Let ${\bf A}$ and ${\bf B}$ be two $n$-dimensional algebras satisfying identities from $T$ and $\mu,\lambda \in \mathbb{L}(T)$ represent ${\bf A}$ and ${\bf B}$ respectively.
We say that ${\bf A}$ {\it degenerates to} ${\bf B}$ and write ${\bf A}\to {\bf B}$ if $\lambda\in\overline{O(\mu)}$.
Note that in this case we have $\overline{O(\lambda)}\subset\overline{O(\mu)}$. Hence, the definition of a degeneration does not depend on the choice of $\mu$ and $\lambda$. If ${\bf A}\not\cong {\bf B}$, then the assertion ${\bf A}\to {\bf B}$ 
is called a {\it proper degeneration}. We write ${\bf A}\not\to {\bf B}$ if $\lambda\not\in\overline{O(\mu)}$.

Let ${\bf A}$ be represented by $\mu\in\mathbb{L}(T)$. Then  ${\bf A}$ is  {\it rigid} in $\mathbb{L}(T)$ if $O(\mu)$ is an open subset of $\mathbb{L}(T)$.
Recall that a subset of a variety is called {\it irreducible} if it cannot be represented as a union of two non-trivial closed subsets. A maximal irreducible closed subset of a variety is called an {\it irreducible component}.
It is well known that any affine variety can be represented as a finite union of its irreducible components in a unique way.
The algebra ${\bf A}$ is rigid in $\mathbb{L}(T)$ if and only if $\overline{O(\mu)}$ is an irreducible component of $\mathbb{L}(T)$.

In the present work we use the methods applied to Lie algebras in \cite{BC99,GRH,GRH2,S90}.
First of all, if ${\bf A}\to {\bf B}$ and ${\bf A}\not\cong {\bf B}$, then $\dim \mathfrak{Der}({\bf A})<\dim \mathfrak{Der}({\bf B})$, where $\mathfrak{Der}({\bf A})$ is the Lie algebra of derivations of ${\bf A}$. We will compute the dimensions of algebras of derivations and will check the assertion ${\bf A}\to {\bf B}$ only for such ${\bf A}$ and ${\bf B}$ that $\dim \mathfrak{Der}({\bf A})<\dim \mathfrak{Der}({\bf B})$. Secondly, if ${\bf A}\to {\bf C}$ and ${\bf C}\to {\bf B}$ then ${\bf A}\to{\bf  B}$. If there is no ${\bf C}$ such that ${\bf A}\to {\bf C}$ and ${\bf C}\to {\bf B}$ are proper degenerations, then the assertion ${\bf A}\to {\bf B}$ is called a {\it primary degeneration}.
It is easy to see that any algebra degenerates to the algebra with zero multiplication. From now on we use this fact without mentioning it.

To prove primary degenerations, we will construct families of matrices parametrized by $t$. Namely, let ${\bf A}$ and ${\bf B}$ be two algebras represented by the structures $\mu$ and $\lambda$ from $\mathbb{L}(T)$ respectively. Let $e_1,\dots, e_n$ be a basis of ${\bf V}$ and $c_{i,j}^k$ ($1\le i,j,k\le n$) be the structure constants of $\lambda$ in this basis. If there exist $a_i^j(t)\in\mathbb{C}$ ($1\le i,j\le n$, $t\in\mathbb{C}^*$) such that $E_i^t=\sum_{j=1}^na_i^j(t)e_j$ ($1\le i\le n$) form a basis of ${\bf V}$ for any $t\in\mathbb{C}^*$, and the structure constants $c_{i,j}^k(t)$ of $\mu$ in the basis $E_1^t,\dots, E_n^t$ satisfy $\lim\limits_{t\to 0}c_{i,j}^k(t)=c_{i,j}^k$, then ${\bf A}\to {\bf B}$. In this case  $E_1^t,\dots, E_n^t$ is called a {\it parametric basis} for ${\bf A}\to {\bf B}$.

If the number of orbits under the action of ${\rm GL}({\bf V})$ on  $\mathbb{L}(T)$ is finite, then the graph of primary degenerations gives the whole picture. In particular, the description of rigid algebras and irreducible components can be easily obtained.

\subsection{The geometric classification of dual mock-Lie algebras}

\subsubsection{Degenerations of $7$-dimensional dual mock-Lie algebras}

\begin{theorem}\label{geotort}
The variety of complex $7$-dimensional dual mock-Lie algebras has
three irreducible component
defined by rigid algebras $\D^7_{09}, \D^7_{13}$ and $\D^7_{14}.$
The complete graph of degenerations in the given variety presented below

\end{theorem}

\bigskip

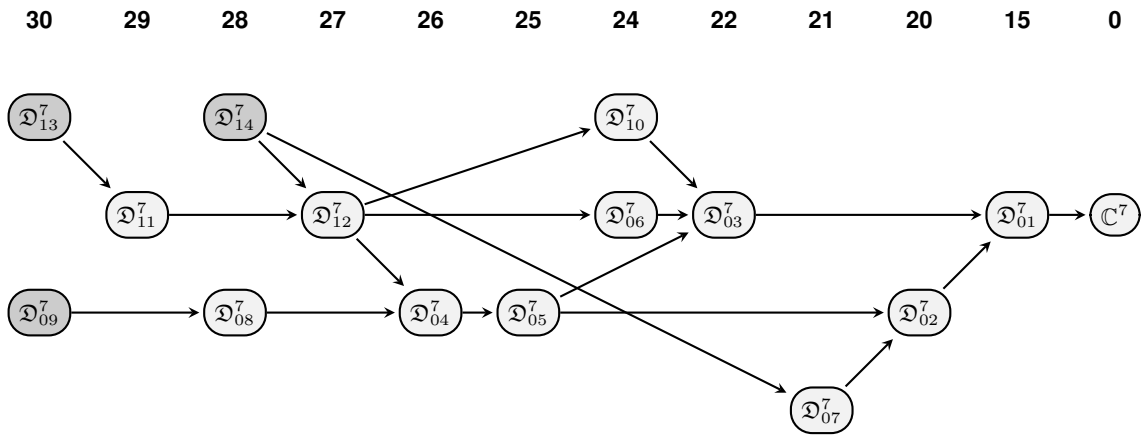
\begin{center}
	\begin{tikzpicture}[->,>=stealth,shorten >=0.05cm,auto,node distance=1.3cm,
	thick,main node/.style={rectangle,draw,fill=gray!10,rounded corners=1.5ex,font=\sffamily \scriptsize \bfseries },rigid node/.style={rectangle,draw,fill=black!20,rounded corners=1.5ex,font=\sffamily \scriptsize \bfseries },style={draw,font=\sffamily \scriptsize \bfseries }]
                    
\node (120)   {30};

\node (110) [right         of=120]      {29};
\node (100) [right         of=110]      {28};
\node (90) [right         of=100]      {27};
\node (80) [right         of=90]      {26};
\node (70) [right         of=80]      {25};
\node (60) [right         of=70]      {24};
\node (50) [right         of=60]      {22};
\node (40) [right         of=50]      {21};
\node (30) [right         of=40]      {20};
\node (20) [right         of=30]      {15};
\node (10) [right         of=20]      {0};
\node (00) [right         of=10]      {};

\node (121)  [below of =120]                      {};
\node (122)  [below of =121]                      { };
\node (123)  [below of =122]                      {};
\node (111)  [below of =110]                      {};
\node (112)  [below of =111]                      { };
\node (113)  [below of =112]                      {};
\node (101)  [below of =100]                      {};
\node (91)  [below of =90]                      {};
\node (92)  [below of =91]                      { };
\node (81)  [below of =80]                      {};
\node (82)  [below of =81]                      { };
\node (71)  [below of =70]                      {};
\node (72)  [below of =71]                      { };
\node (73)  [below of =72]                      {};
\node (61)  [below of =60]                      {};
\node (62)  [below of =61]                      { };
\node (51)  [below of =50]                      {};
\node (52)  [below of =51]                      { };
\node (53)  [below of =52]                      { };
\node (41)  [below of =40]                      {};
\node (42)  [below of =41]                      { };
\node (31)  [below of =30]                      {};
\node (32)  [below of =31]                      { };
\node (21)  [below of =20]                      {};
\node (22)  [below of =21]                      { };
\node (11)  [below of =10]                      {};
\node (12)  [below of =11]                      { };

\node [rigid node] (27b)  [below of =122]                      {$\D^7_{09}$};

\node [rigid node] (37d)  [below of =120]                      {$\D^7_{13}$};
\node (yy)  [right of =37d]                      { };
\node [rigid node] (d14)  [right of =yy]                      {$\D^7_{14}$};
\node [main node] (37b)  [below of =111]                      {$\D^7_{11}$};
\node (xx)  [right of =27b]                      { };
\node [main node] (27a)  [right of =xx]                      {$\D^7_{08}$};
\node [main node] (37c)  [below of =91]                      {$\D^7_{12}$};
\node [main node] (n31n31)  [below of =82]                      {$\D^7_{04}$};
\node [main node] (37a)  [below of =60]                      {$\D^7_{10}$};
\node [main node] (n62)  [below of =72]                      {$\D^7_{05}$};
\node [main node] (n61)  [below of =61]                      {$\D^7_{06}$};

\node  (41)  [below of =40]                      { };
\node  (42)  [below of =41]                      { };
\node  (43)  [below of =42]                      { };
\node [main node] (17)  [below of =43]                      {$\D^7_{07}$};
\node [main node] (n51)  [below of =51]                      {$\D^7_{03}$};

\node [main node] (n53)  [below of =32]                      {$\D^7_{02}$};

\node [main node] (n31)  [below of =21]                      {$\D^7_{01}$};

\node [main node] (c)  [below of =11]                      {$\mathbb{C}^7$};

\path[every node/.style={font=\sffamily\tiny}]

(n31)  edge [bend right=0, color=black] node{}  (c)

(d14)  edge [bend right=0, color=black] node{}  (37c)

(n53)  edge [bend right=0, color=black] node{}  (n31)

(n51)  edge [bend right=0, color=black] node{}  (n31)

(17)  edge [bend right=0, color=black] node{}  (n53)

(n61)  edge [bend right=0, color=black] node{}  (n51)

(37a)  edge [bend right=0, color=black] node{}  (n51)

(n62)  edge [bend right=0, color=black] node{}  (n51)
(n62)  edge [bend right=0, color=black] node{}  (n53)

(n31n31)  edge [bend right=0, color=black] node{}  (n62)

(37c)  edge [bend right=0, color=black] node{}  (37a)

(37c)  edge [bend right=0, color=black] node{}  (n61)
(37c)  edge [bend right=0, color=black] node{}  (n31n31)

(27a)  edge [bend right=0, color=black] node{}  (n31n31)

(d14)  edge [bend right=0, color=black] node{}  (17)

(37b)  edge [bend right=0, color=black] node{}  (37c)

(27b)  edge [bend right=0, color=black] node{}  (27a)

(37d)  edge [bend right=0, color=black] node{}  (37b)
;
\end{tikzpicture}
\end{center}

\begin{Proof}
Thanks to \cite{ale2} we have all degenerations in the variety of all $7$-dimensional $2$-step nilpotent  Lie algebras. By some easy calculation, we have that the dimension of the algebra of derivation of the algebra $\D$
is 21.
Hence, it can not degenerates to $\D^7_{08}, \D^7_{09}, \D^7_{11}, \D^7_{13}.$

The degeneration $\D^7_{14} \to \D^7_{12}$ is obtained by the following parametric basis  

$$\begin{array}{llll}
E^t_1= t e_4 &  E^t_2= t^2 e_2 - e_3 & \multicolumn{2}{l}{E^t_3= t e_3 + t e_5 + t^3 e_6}\\
E^t_4= e_1 + e_2 + t^2 e_4 - e_5 & E^t_5= t e_7 & E^t_6= t^3 e_6  & E^t_7= e_5 +  e_6.\\
\end{array}$$

The degeneration $\D^7_{14} \to \D^7_{07}$ is obtained by the following parametric basis

$$\begin{array}{llll}

E^t_1= t e_1 &
E^t_2=  e_6 &
E^t_3= e_2 &
E^t_4= -t e_5 \\
E^t_5= t e_3 &
E^t_6= e_4 &
E^t_7= t e_7

\end{array}$$

\end{Proof}

\begin{remark}
Note that, the graph of primary degenerations of  $7$-dimensional $2$-step nilpotent Lie algebras from  \cite{ale2} is not correct. 
We gave the corrected graph of degenerations of $7$-dimensional $2$-step nilpotent Lie algebras from  \cite{ale2}. 
\end{remark}

\subsubsection{The geometric classification  of $8$-dimensional dual mock-Lie algebras}
Thanks to \cite{ale} we have that the variety of $8$-dimensional $2$-step nilpotent Lie algebras has 
three rigid algebras: $\D_{17}^8, \D_{30}^8$ and $\D_{33}^8.$
It is easy to see, that the algebra $\D_{36}^8$ is satisfying 
the following invariant conditions $A_4A_5=0$ and $A_1A_4 \subseteq A_8,$
but the cited algebras are not satisfy it.
It is follow that there are no degenerations $\D_{36}^8 \to  \D_{17}^8, \D_{30}^8, \D_{33}^8.$

The degeneration $\D^8_{36} \to \D^8_{14}$ is obtained by the following parametric basis  
$$\begin{array}{llllllll}
E^t_1=  e_1 &
E^t_2=  e_2 &
E^t_3=  e_3 &
E^t_4=  e_4 &
E^t_5=  e_5 &
E^t_6=  e_6 &
E^t_7= e_8 & 
E^t_8= t e_7. 
\end{array}$$
Hence, we have the following theorem

\begin{theorem}\label{geotort}
The variety of complex $8$-dimensional dual mock-Lie algebras has
four  irreducible component
defined by rigid algebras $\D^8_{17}, \D^8_{30}, \D^8_{33}$ and $\D^8_{36}.$

\end{theorem}


\begin{thebibliography}{99}

\bibitem{ack}
Abdelwahab H.,  Calder\'on A.J., Kaygorodov I.,
    The algebraic and geometric classification of nilpotent binary Lie algebras, 
      International Journal of Algebra and Computation,  29 (2019), 6, 1113--1129
  
    
\bibitem{omirov}
Adashev J.,  Camacho L.,  Omirov B.,
    Central extensions of null-filiform and naturally graded filiform non-Lie Leibniz algebras,
    Journal of Algebra, 479 (2017), 461--486.
    
\bibitem{am} Agore A.L., Militaru G. On a type of commutative algebras, 
    Linear Algebra and its Applications,  485 (2015), 222--249.

\bibitem{ale}
Alvarez M.A., 
    On rigid $2$-step nilpotent Lie algebras, 
    Algebra Colloquium,  25  (2018), 2, 349--360.


\bibitem{ale2}
Alvarez M.A., 
    The variety of $7$-dimensional $2$-step nilpotent Lie algebras, 
    Symmetry,  10 (2018), 1, 26.
    
    \bibitem{maria}
Alvarez M.A., Hern\'{a}ndez I., Kaygorodov I.,
    Degenerations of Jordan superalgebras,
    Bulletin of the Malaysian Mathematical Sciences Society,  42 (2019),   6, 3289--3301.



\bibitem{bf}
 Burde D., Fialowski A, Jacobi–Jordan algebras, 
 Linear Algebra and its Applications,  459 (2014), 586--594.

\bibitem{BC99} 
Burde D., Steinhoff C.,
    Classification of orbit closures of $4$--dimensional complex Lie algebras,
    Journal of Algebra, 214 (1999), 2, 729--739.

\bibitem{degr3}
Cicalò S., De Graaf W.,   Schneider C.,
    Six-dimensional nilpotent Lie algebras,
    Linear Algebra and its Applications, 436 (2012), 1, 163--189.

\bibitem{usefi1}
Darijani I., Usefi H.,
    The classification of 5-dimensional $p$-nilpotent restricted Lie algebras over perfect fields. I,
    Journal of Algebra, 464 (2016), 97--140.


\bibitem{degr2}
De Graaf W., 
    Classification of 6-dimensional nilpotent Lie algebras over fields of characteristic not $2$, 
    Journal of Algebra, 309  (2007), 2, 640--653.

\bibitem{degr1}
De Graaf W., 
    Classification of nilpotent associative algebras of small dimension,
    International Journal of Algebra and Computation, 28 (2018),  1, 133--161.

\bibitem{fkkv}
Fern\'andez Ouaridi A.,  Kaygorodov I.,  Khrypchenko M., Volkov Yu., 
    Degenerations of nilpotent algebras,
    arXiv:1905.05361

\bibitem{gk}
Getzler E., Kapranov M., 
    Cyclic operads and cyclic homology, 
    Geometry, Topology and Physics for Raoul Bott (ed. S.-T. Yau), International Press, 1995, 167--201.

\bibitem{gkk}
Gorshkov I., Kaygorodov I., Khrypchenko M.,
    The algebraic classification of nilpotent Tortkara algebras,
    arXiv:1904.00845.

\bibitem{gkk19}
Gorshkov I., Kaygorodov I., Khrypchenko M.,
    The geometric classification of nilpotent Tortkara algebras,
    Communications in Algebra, 2019, DOI:	10.1080/00927872.2019.1635612


\bibitem{gkks}
Gorshkov I., Kaygorodov I., Kytmanov A., Salim M.,
    The variety of nilpotent Tortkara algebras,
    Journal of Siberian Federal University. Mathematics \& Physics, 12 (2019), 2, 173--184.

\bibitem{GRH}
Grunewald F.,  O'Halloran J.,
    Varieties of nilpotent Lie algebras of dimension less than six,
    Journal of Algebra, 112 (1988), 2, 315--325.

\bibitem{GRH2}
Grunewald F., O'Halloran J.,
    A Characterization of orbit closure and applications,
    Journal of Algebra, 116 (1988), 1, 163--175.


\bibitem{ha16}
Hegazi A., Abdelwahab H.,
    Classification of five-dimensional nilpotent Jordan algebras,
    Linear Algebra and its Applications, 494 (2016), 165--218.



\bibitem{hac16}
Hegazi A., Abdelwahab H., Calderón Martín A.,
    The classification of $n$-dimensional non-Lie Malcev algebras with $(n-4)$-dimensional annihilator, 
    Linear Algebra and its Applications, 505 (2016), 32--56.

\bibitem{hac18}
Hegazi A., Abdelwahab H.,  Calderón Martín A.,
    Classification of nilpotent Malcev algebras of small dimensions over arbitrary fields of characteristic not $2$,
    Algebras and  Represention Theory, 21 (2018), 1, 19--45.

\bibitem{hjs} Hentzel I., Jacobs D., Sverchkov S., 
    On exceptional nil of index 3 Jordan algebras, 
    Preprint, Novosibirsk State University, 1997.


\bibitem{ikv17}
Ismailov N., Kaygorodov I.,  Volkov Yu.,
    The geometric classification of Leibniz algebras,
    International Journal of Mathematics, 29  (2018), 5, 1850035.

\bibitem{ikv18}
Ismailov N., Kaygorodov I.,  Volkov Yu.,
    Degenerations of Leibniz and anticommutative algebras,
    Canadian Mathematical Bulletin, 62 (2019),  3, 539--549.



\bibitem{kkk18}
Karimjanov I., Kaygorodov I., Khudoyberdiyev A.,
    The algebraic and geometric classification of nilpotent Novikov  algebras, 
    Journal of Geometry and Physics, 143 (2019), 11--21.



\bibitem{kkl19}
Kaygorodov I., Khrypchenko M., Lopes S.,
    The algebraic and geometric classification of nilpotent anticommutative  algebras, 
    preprint



\bibitem{kpv}
Kaygorodov I., Popov Yu., Volkov Yu.,
    Degenerations of binary-Lie and nilpotent Malcev algebras,
    Communications in Algebra, 46 (2018), 11, 4929--4941.


\bibitem{kv16}
Kaygorodov I.,   Volkov Yu.,
    The variety of $2$-dimensional algebras over an algebraically closed field,
    Canadian  Journal of Mathematics,  71 (2019),  4, 819--842.



\bibitem{kv17}
Kaygorodov I., Volkov Yu., 
    Complete classification of algebras of level two,  
    Moscow Mathematical Journal, 19 (2019), 3,  485--521.
    

\bibitem{mr} Markl M., Remm E., (Non-)Koszulness of operads for n-ary algebras, galgalim and other curiosities, 
    Journal of Homotopy and Related Structures, 10 (2015), 939--969.

\bibitem{ok} Okubo S., Kamiya N, Jordan--Lie super algebra and Jordan-Lie triple system, 
    Journal of Algebra, 198 (1997), 388--411.


\bibitem{ren}
Ren B.,  Zhu L., 
    Classification of $2$-step nilpotent Lie algebras of dimension $9$ with $2$-dimensional center,   
    Czechoslovak Mathematical Journal,  67 (2017), 4, 953--965.

\bibitem{S90}
Seeley C.,  
    Degenerations of 6-dimensional nilpotent Lie algebras over $\mathbb{C}$, 
    Communications in Algebra, 18 (1990), 3493--3505.


\bibitem{ss78}
Skjelbred T., Sund T.,
    Sur la classification des algebres de Lie nilpotentes,
    C. R. Acad. Sci. Paris Ser. A-B, 286 (1978), 5,  A241--A242.




\bibitem{w} Walcher S., 
    On algebras of rank three, 
    Communications in Algebra, 27 (1999), 3401--3438

\bibitem{zhev}
    Zhevlakov K., 
    Solvability and nilpotence of Jordan rings, Algebra i Logika, 5 (1966),  37--58 (in Russian).

\bibitem{zsss} Zhevlakov K., Slin’ko A.M., Shestakov I.P.,  Shirshov A.I., Rings That Are Nearly Associative, Nauka, Moscow,
1978 (in Russian); Academic Press, 1982 (English translation).

\bibitem{zusmanovich}
Zusmanovich P., 
    Central extensions of current algebras,
    Transactions of the American Mathematical Society, 334 (1992),  1, 143--152.


\bibitem{pasha}
Zusmanovich P., 
    Special and exceptional mock-Lie algebras,
    Linear Algebra and its Applications,  518 (2017), 79--96. 

\end{thebibliography}
\end{document}